\documentclass[11pt,leqeq]{article}
\usepackage{amssymb,amsthm}
\usepackage{amsmath}
\usepackage{amscd}
\def\ob{{\rm{Ob}}}

\newtheorem{thm}{Theorem}
\newtheorem{lema}[thm]{Lemma}

\newtheorem{defi}[thm]{Definition}

\newtheorem{problem}[thm]{Problem}

\begin{document}
\title{On the combinatorics of hypergeometric functions}
\author{H\'{e}ctor Bland\'{i}n \
and \ Rafael D\'{i}az} \maketitle
\date

\begin{abstract}
\noindent  We give combinatorial interpretation for
hypergeometric functions associated with tuples of rational
numbers.\\

\noindent AMS Subject Classification: \ \ 05A99, 33C20, 18A99.\\
\noindent Keywords:\ \ Rational Combinatorics, Groupoids, Generalized
Hypergeometric Functions.

\end{abstract}
\section{Introduction}
Given a natural number $n \in \mathbb{N}$ and a real number $a
\in \mathbb{R}$ the Pochhammer symbol $(a)_n$ is defined by the
identity $(a)_n = a(a + 1)(a + 2) . . .(a + n - 1)$. Later we
shall also need a generalized form of the Pochhammer symbol,
introduced in \cite{DP}, defined as follows: let $k
\in
\mathbb{R}$ then $$ (a)_{n,k} = a(a + k)(a + 2k)...(a + (n -
1)k).$$

Now given a tuple of complex numbers $a_{1},...,a_{k}$ and another
tuple of complex numbers $b_{1},...,b_{l}$ such that no $b_i$ is a
negative integer or zero, then there is an associated
hypergeometric series  $h(a_{1},...,a_{k};b_{1},...,b_{l}) \in
\mathbb{C}[[x]],$ the formal power series with complex
coefficients defined by:
$$h(a_{1},...,a_{k};b_{1},...,b_{l}) = \sum_{n=0}^{\infty}
\frac{(a_{1})_{n}...(a_{k})_{n}}{(b_{1})_{n}...(b_{l})_{n}} \frac{x^{n}}{n!}.$$
\

The goal of this paper is to find a combinatorial interpretation
of the coefficients
$\frac{(a_{1})_{n}...(a_{k})_{n}}{(b_{1})_{n}...(b_{l})_{n}}$ of
the hypergeometric series $h(a_{1},...,a_{k};b_{1},...,b_{l})$.
Without further restrictions this problem seems hopeless since in
general the coefficients
$\frac{(a_{1})_{n}...(a_{k})_{n}}{(b_{1})_{n}...(b_{l})_{n}}$ are
complex numbers and, to this day, no good definition of
"combinatorial interpretation" for complex numbers is known. \\

If we demand that $a_{1},...,a_{k}$ and $b_{1},...,b_{l}$ be
rational numbers, then clearly the coefficients
$\frac{(a_{1})_{n}...(a_{k})_{n}}{(b_{1})_{n}...(b_{l})_{n}}$ will
also be rational numbers. In \cite{BD} we proposed a general
setting in which the phrase "combinatorial interpretation of a
rational number" is well defined and makes good sense. Indeed a
combinatorial interpretation of Bernoulli and Euler numbers is
provided in \cite{BD}. So working within that setting the question
of the combinatorial interpretation of the coefficients of
hypergeometric functions associated with
tuples of rational numbers makes sense. \\

\noindent The rest of this paper is devoted to, first, explaining what
do we mean when we talk about  combinatorial interpretations of
rational numbers, and second, use this notion to uncover the
combinatorial meaning of the hypergeometric series
$h(a_{1},...,a_{k};b_{1},...,b_{l})$ with
$a_{1},...,a_{k},b_{1},...,b_{l} \in \mathbb{Q^+}$. Notice that we
are restricting our attention to positive rational numbers. We do
this for simplicity, the case of negative rational numbers can
also be dealt with  our methods, and will be the subject of future
work.

\section{Rational Combinatorics}

In this section we explain what we mean when we talk about the
combinatorial interpretation of a formal power series in
$\mathbb{Q^+}[[x]]$. We only give the material needed to attack
the  problem of this paper, namely, the combinatorial
interpretation of $h(a_{1},...,a_{k};b_{1},...,b_{l})$ with
$a_{1},...,a_{k},b_{1},...,b_{l} \in \mathbb{Q^+}$. The interested
reader may find more information in references \cite{BaezDolan}
\cite{Bergeron}, \cite{BD} and \cite{RDEP}.\\

Let us define what we mean by the combinatorial interpretation of
a positive rational number. To do so we start by introducing the
category $gpd$ of finite groupoids. Recall that a category $C$
consists of a collection of objects $Ob(C)$, and a set of
morphisms $C(x,y)$ for each pair of objects $x,y \in C$, together
with a natural list of axioms  \cite{SMacLane}. Let us recall the
definition of a finite groupoid
\cite{BaezDolan}.

\begin{defi}
A groupoid $G$ is a category such that all its morphisms are
invertible. A groupoid $G$ is called finite if $Ob(G)$ is a finite
set and $G(x,y)$ is a finite set for all $x,y\in Ob(G).$
\end{defi}

\begin{defi}
The category $gpd$  is such that its objects $Ob(gpd)$ are finite
groupoids. For groupoids $G,H \in Ob(gpd)$ morphisms in $gpd(G,H)$
are functors $F:G \rightarrow H.$
\end{defi}

There are two rather simple examples of a finite groupoid: 1) Any
finite set $x$ may be regarded as a groupoid. The objects of $x$
is $x$ itself and the only morphisms in $x$ are the identities. 2)
Any finite group $G$ may be regarded as a groupoid $\overline{G}$.
$\overline{G}$ has only one object $1$ and
$\overline{G}(1,1)=G.$\\

The category $gpd$ has several remarkable properties which we
briefly summarize:

\begin{itemize}
\item There is a bifunctor $\bigsqcup: gpd \times gpd
\longrightarrow gpd$ called disjoint union and defined for $G,H$
in $gpd$ by $Ob(G \sqcup H) = Ob(G) \sqcup Ob(H)$ and for $x,y\in
Ob(G \sqcup H)$ the morphisms from $x$ to $y$ are given by
$$G \sqcup H(x,y)=\left\{\begin{array}{cc}
G(x,y) & \mathrm{if}\ x,y\in Ob(G),\\
H(x,y) & \mathrm{if}\ x,y\in Ob(H),\\
\emptyset   & otherwise.
\end{array}\right.$$

\item There is a bifunctor $\times:gpd \times gpd \longrightarrow gpd$
called Cartesian product and defined for $G,H$ in $gpd$ by $Ob(G
\times H) = Ob(G) \times Ob(H)$ and for all
$(x_{1},y_{1}),(x_{2},y_{2})\in Ob(G \times H)$ we have $G
\times H((x_{1},y_{1}),(x_{2},y_{2}))=G(x_{1},x_{2})\times H(y_{1},y_{2}).$

\item There is an empty groupoid $\emptyset$. It is such that
$Ob(\emptyset) = \emptyset$. $\emptyset$ is the neutral element in
$gpd$ with respect to disjoint union.

\item There is a groupoid $1$ such that it has only one object
and only one morphism, the identity, between that object and
itself. $1$ is an identity with respect to Cartesian product in
$gpd$.

\item There is a valuation map $|\mbox{  } |: Ob(gpd) \rightarrow
\mathbb{Q^{+}}$, given on $G \in Ob(gpd)$ by
$$|G|=\sum_{x \in D(G)}\frac{1}{|G(x,x)|},$$ where $D(G)$ denotes the set
of isomorphisms classes of objects of $G$ and for $|A|$ denotes
the cardinality of $A.$ $|G|$ is called the cardinality of the
groupoid $G$.

\item The valuation map $|\mbox{  } |: Ob(gpd) \rightarrow
\mathbb{Q^{+}}$ is such that $|G \sqcup H| = |G| + |H|$, $|G \times
H|=|G||H|$ for all $G,H \in gpd$. Also $|\emptyset|=0$ and
$|1|=1$.
\end{itemize}

The properties above suggest the following

\begin{defi}
If $a=|G|$, where $a \in \mathbb{Q}^{+}$ and $G \in gpd$, then we
say that $a$ is the cardinality of the groupoid $G$, and also that
$G$ is a combinatorial interpretation (in terms of finite
groupoids) of the rational number $a$.
\end{defi}

Notice that there are many groupoids with the same cardinality, so
the combinatorial interpretation of a rational number is by no
means unique. This should not be seen as something negative,
indeed the freedom to choose among several different possible
interpretations is what makes the subject of enumerative
combinatorics ( either for integers or for rationals ) so
fascinating. For example the reader may look at \cite{Stanley}
where a list of references with over $30$ different combinatorial
interpretations of the Catalan numbers is given.\\

For example, in this paper we choose the combinatorial
interpretation of the number $\frac{1}{n}$ to be
$\overline{\mathbb{Z}_n}$ the groupoid associated to the cyclic
group with $n$-elements. We believed our choice is justified for
simplicity. Nevertheless, this seemingly naive choice is actually
quite subtle. For example, we think that the number $\frac{1}{nm}$
is better interpreted by $\overline{\mathbb{Z}}_n
\times
\overline{\mathbb{Z}}_m$ than by $\overline{\mathbb{Z}}_{nm}$. Since
$\overline{\mathbb{Z}}_n \times
\overline{\mathbb{Z}}_m$ and $\overline{\mathbb{Z}}_{nm}$ are not isomorphic groupoids they
provide really different combinatorial interpretations for
$\frac{1}{nm}$. Notice that in general we have that
$|\overline{\mathbb{Z}}_{b}^{\oplus a}|=\frac{a}{b}.$\\

Now that we have defined what we mean by a combinatorial
interpretation (in terms of finite groupoids) of a rational
number, we proceed to define what is a combinatorial
interpretation of a formal power series with positive rational
coefficients. First we introduce the groupoid $\mathbb{B}$.
Objects of $\mathbb{B}$ are finite sets. For finite sets $x$ and
$y$, morphisms in $\mathbb{B}(x,y)$ are bijections $f: x
\rightarrow y$.

\begin{defi} The category of
$\mathbb{B}$-$gpd$ species is the category ${gpd}^{\mathbb{B}}$ of
functors from $\mathbb{B}$ to $gpd$. An object of
${gpd}^{\mathbb{B}}$  is called a $\mathbb{B}$-$gpd$ species or a
species.
\end{defi}

The category of rational species ${gpd}^{\mathbb{B}}$ has the
following properties

\begin{itemize}
\item There is a bifunctor $+:{gpd}^{\mathbb{B}} \times {gpd}^{\mathbb{B}}
\longrightarrow {gpd}^{\mathbb{B}}$ given for $F,G \in
{gpd}^{\mathbb{B}}$ as follows: $F+G(x)=F(x) \sqcup G(x)$ for all
$x
\in Ob(\mathbb{B})$.  Bifunctor $+$ is called the sum of species.

\item There is a bifunctor $\times:{gpd}^{\mathbb{B}} \times {gpd}^{\mathbb{B}}
\longrightarrow {gpd}^{\mathbb{B}}$ given for $F,G \in
{gpd}^{\mathbb{B}}$ as follows $F \times G(x) = F(x) \times G(x)$
for all $x \in Ob(\mathbb{B})$. The bifunctor $\times$ is called
the Hadamard product.

\item There is a bifunctor ${gpd}^{\mathbb{B}} \times {gpd}^{\mathbb{B}}
\longrightarrow {gpd}^{\mathbb{B}}$  given for $F,G \in
{gpd}^{\mathbb{B}}$ as follows: $FG(x)= \bigsqcup_{a \subseteq
x}F(a)
\times G(x\backslash a)$ for all $x
\in Ob(\mathbb{B})$.  This bifunctor is called the product of species.

\item Let ${gpd}_{0}^{\mathbb{B}}$ be the full subcategory of species
$F$ such that $F( \emptyset)= \emptyset$. There is a bifunctor
$$\circ :{gpd}^{\mathbb{B}} \times {gpd}_{0}^{\mathbb{B}}
\longrightarrow {gpd}^{\mathbb{B}}$$ given for $F \in
{gpd}^{\mathbb{B}}$ and $G \in {gpd}_{0}^{\mathbb{B}}$ by $$F
\circ G(x)= \bigsqcup_{p \in
\Pi(x)} F(p)\prod_{b \in p}G(b)$$  for all $x
\in Ob(\mathbb{B})$. Above $\Pi(x)$ denotes the set of partitions of $x$.
This bifunctor is called the composition or substitution of
species.

\item The species $1$ is defined by $1(x)= 1$ if $|x|=1$ and
$1(x)=\emptyset$ otherwise. The species $0$ is defined by $0(x)=
\emptyset$ if $x \neq \emptyset$ and $0(\emptyset)=1$.

\item There is a valuation map $|\mbox{  } |: Ob(gpd^{\mathbb{B}}) \rightarrow
\mathbb{Q^{+}}[[x]]$ given on $F \in Ob(gpd^{\mathbb{B}})$ by
$$|F|=\sum_{n \in \mathbb{N}}|F[n]| \frac{x^n}{n!},$$ where
$[n]=\{ 1,2,...,n \}$ and $[0]=\emptyset$.

\item The valuation map $|\mbox{  } |: Ob(gpd^{\mathbb{B}}) \rightarrow
\mathbb{Q^{+}}[[x]]$ satisfies the following properties: $|F + G| =
|F| + |G|, |F \times G|=|F| \times |G|,$ the Hadamard product of
powers series, i.e., coefficientwise multiplication, $|FG|=|F||G|,
|F \circ G|= |F|(|G|), |1|=1,$ and $|0|=0.$
\end{itemize}

The properties above allow one to make the following

\begin{defi}
If $f=|F|$, where $f \in \mathbb{Q^{+}}[[x]]$ and $F \in
Ob(gpd^{\mathbb{B}}),$ then we say that $f$ is the generating
series associated with the species $F$, and also that the species
$F$ is a combinatorial interpretation (in terms of functors from
finite sets to groupoids) of the formal power series $f$.
\end{defi}

\smallskip

For example let $\mathbb{Z}:\mathbb{B}\rightarrow gpd$ be the
species sending $x\in\ob(\mathbb{B})$  into the groupoid
$\mathbb{Z}(x)$ such that
$\ob(\mathbb{Z}(x))=\left\{\begin{array}{cc} \{x\} & if\ x\neq\phi,\\
\phi & if\ x=\phi.\end{array}\right.$ and
$\mathbb{Z}(x)(x,x)=\mathbb{Z}_{|x|}$. Then
\begin{equation*}
|\mathbb{Z}|=\sum_{n=0}^{\infty}|\mathbb{Z}(n)|\frac{x^{n}}{n!}=
\sum_{n=1}^{\infty}\frac{1}{n}\frac{x^{n}}{n!}=\int\frac{e^{x}-1}{x}dx.
\end{equation*}

\smallskip

After we have made all these remarks we can rephrase our original
problem of finding a combinatorial interpretation for the
hypergeometric series as follows: given tuples of positive
rational numbers $a_1,...,a_k,b_1,...,b_l,$ find a  species
$H(a_1,...,a_k;b_1,...,b_l):\mathbb{B} \rightarrow gpd$ such that
$|H(a_1,...,a_k;b_1,...,b_l)|=h(a_1,...,a_k;b_1,...,b_l)$. A
solution to this problem will be given in the next section. We
remark that at this point our extension of the notion of
combinatorial species introduced by Joyal in \cite{j1} and
\cite{j2}, replacing the target category from $\mathbb{B}$ to
$gpd$ becomes necessary. The generating series associated with a
combinatorial species is a formal power series with integer
coefficients, and thus in that formalisms there is no hope to find
functors $H(a_1,...,a_k;b_1,...,b_l)$ with the desired properties.

\section{Combinatorial interpretation of hypergeometric functions}

In this section we construct functors
$$H(\frac{a_1}{b_1},...,\frac{a_k}{b_k};\frac{c_1}{d_1},...,\frac{c_l}{d_l}):\mathbb{B} \rightarrow gpd$$
for $a_1,...,a_k, b_1,...,b_k,c_1,...,c_l,d_1,...,d_l
\in \mathbb{N}^{+}$ such that
$$|H(\frac{a_1}{b_1},...,\frac{a_k}{b_k};\frac{c_1}{d_1},...,\frac{c_l}{d_l})|=
h(\frac{a_1}{b_1},...,\frac{a_k}{b_k};\frac{c_1}{d_1},...,\frac{c_l}{d_l}).$$
The construction of functors
$H(\frac{a_1}{b_1},...,\frac{a_k}{b_k};\frac{c_1}{d_1},...,\frac{c_l}{d_l})$
provides the combinatorial interpretation for the hypergeometric
functions
$h(\frac{a_1}{b_1},...,\frac{a_k}{b_k};\frac{c_1}{d_1},...,\frac{c_l}{d_l})$
promised in the introduction. As remarked above combinatorial
interpretations are by no means unique. Nevertheless they are not
arbitrary, and a good combinatorial interpretation is subject to
many constrains. In our construction of
$H(\frac{a_1}{b_1},...,\frac{a_k}{b_k};\frac{c_1}{d_1},...,\frac{c_l}{d_l})$
we are going to carefully respect the Hadamard product structure
implicit in the definition of
$h(\frac{a_1}{b_1},...,\frac{a_k}{b_k};\frac{c_1}{d_1},...,\frac{c_l}{d_l}).$
Namely, since
$$h(\frac{a_1}{b_1},...,\frac{a_k}{b_k};\frac{c_1}{d_1},...,\frac{c_l}{d_l})
= \times_{i=1}^{k}h(\frac{a_i}{b_i};\emptyset)\times
\times_{j=1}^{l}h(\emptyset,\frac{c_j}{d_j}),$$
it is natural to demand that

$$H(\frac{a_1}{b_1},...,\frac{a_k}{b_k};\frac{c_1}{d_1},...,\frac{c_l}{d_l})
= \times_{i=1}^{k}H(\frac{a_i}{b_i};\emptyset)\times
\times_{j=1}^{l}H(\emptyset,\frac{c_j}{d_j}),$$
where $|H(\frac{a_i}{b_i};\emptyset)| =
h(\frac{a_i}{b_i};\emptyset)$ and $|H(\emptyset,\frac{c_j}{d_j})|=
h(\emptyset,\frac{c_j}{d_j})$ for all $ 1 \leq i \leq k$ and $1
\leq j \leq l.$ We begin by introducing some definitions and notations.
Next we introduce the analogue in the category $gpd$ of the
generalized Pochhammer symbol.

\begin{defi}\label{pf}
Fix a finite groupoid $K$. The functorial Pochhammer symbol $( \ \
)_{n,K}: gpd \rightarrow gpd$ is given by
$$( G )_{n,K}=\prod_{i=0}^{n-1}(G
\sqcup (K \times [i])),$$ for $G$ in $gpd$. If $K=1$
we write $( \ \ )_{n}$ instead of  $( \ \ )_{n,1}.$
\end{defi}

With this notation we have the following

\begin{lema} \label{gps}
If $|K|=k$ and $|G|=g$, then $|( G )_{n,K}|=(g)_{n,k}$.
\end{lema}

\begin{proof}
$|( G )_{n,K}|=|\prod_{i=0}^{n-1}(G
\sqcup (K \times [i]))|= \prod_{i=0}^{n-1}(|G|
+ |K|i)= (g)_{n,k}.$
\end{proof}

For example one can let the set
$Ob((\overline{\mathbb{Z}}_{b}^{\oplus a})_n)$  be
$[a]\times[a+1]\times ... \times[a+n-1]$. An object $x \in
Ob((\overline{\mathbb{Z}}_{b}^{\oplus a})_n)$ is a tuple $x =
(x_1,...,x_n)$ such that $x_i \in \mathbb{N}$ and $1
\leq x_i \leq a+i-1.$ There are morphisms only from an object to
itself and $$(\overline{\mathbb{Z}}_{b}^{\oplus a})_{n}(x,x)=
\overline{\mathbb{Z}}_{b}^{\times c(x)},$$ where  $c(x)=|\{i \ | \  1
\leq x_i \leq a \}|$ for $x \in Ob((\overline{\mathbb{Z}}_{b}^{\oplus
a})_n).$ Clearly $|(\overline{\mathbb{Z}}_{b}^{\oplus
a})_n|=(\frac{a}{b})_n$.

\begin{defi}
For $a,b \in \mathbb{N}^{+}$ we let the functor
$H(\frac{a}{b},\emptyset):\mathbb{B} \rightarrow gpd$ be such that
on $x \in \mathbb{B}$ it is given by
$$H(\frac{a}{b},\emptyset)(x)=([a])_{|x|,[b]}\overline{\mathbb{Z}}_{b}^{x}.$$
\end{defi}

Explicitly $Ob(H(\frac{a}{b},\emptyset)(x))= [a]\times[a+b] \times
... \times [a +(|x|-1)b],$  and for $s \in
Ob(H(\frac{a}{b},\emptyset)(x))$ we have that
$H(\frac{a}{b},\emptyset)(x)(s,s)=\mathbb{Z}_{b}^{x}.$

\begin{lema}\label{ya}
For $a,b \in \mathbb{N}^{+}$ we have that
$|H(\frac{a}{b},\emptyset)|=h(\frac{a}{b},\emptyset)=
\frac{1}{(1 - x)^{\frac{a}{b}}}.$
\end{lema}

\begin{proof}
For $n \in \mathbb{N}$ one obtains
$|H(\frac{a}{b},\emptyset)([n])|=|([a])_{n,[b]}||\overline{\mathbb{Z}}_{b}^{[n]}|=
\frac{(a)_{n,b}}{b^{n}}=(\frac{a}{b})_{n}.$ This identity implies the
desired result.
\end{proof}

Another combinatorial interpretation for
$h(\frac{a}{b},\emptyset)$ is given by

\begin{lema}
The functor sending a set $x$ into
$(\overline{\mathbb{Z}}_{b}^{\oplus a})_{|x|}$ provides a
combinatorial interpretation for $h(\frac{a}{b},\emptyset)$.
\end{lema}

\begin{defi}
For $c,n,d \in  \mathbb{N}$ we set
$\overline{\mathbb{Z}}_{c,n,d}=\prod_{i=0}^{n-1}\overline{\mathbb{Z}}_{c
+ id}.$
\end{defi}

\begin{lema}
For $c,n,d \in  \mathbb{N}$ we have that
$|\overline{\mathbb{Z}}_{c,n,d}|=
\frac{1}{(c)_{n,d}}$.
\end{lema}

\begin{defi}
For $c,d \in \mathbb{N}^{+}$ we let the functor
$H(\emptyset;\frac{c}{d}):\mathbb{B} \rightarrow gpd$ be such that
on $x
\in \mathbb{B}$ it is given by
$$H(\emptyset;\frac{c}{d})(x)=[d]^{x}\overline{\mathbb{Z}}_{c,|x|,d}.$$
\end{defi}

\begin{lema}
For $c,d \in \mathbb{N}^{+}$ we have that
$|H(\emptyset;\frac{c}{d})|=h(\emptyset;\frac{c}{d}).$
\end{lema}

\begin{proof}
For $n \in \mathbb{N}$ one obtains
$|H(\emptyset;\frac{c}{d})([n])|=|[d]^{[n]}\overline{\mathbb{Z}}_{c,n,d}|=
\frac{d^{n}}{(c)_{n,d}}=(\frac{c}{d})_{n}.$ This identity implies the desired result.
\end{proof}

We are now ready to define the functor
$H(\frac{a_1}{b_1},...,\frac{a_k}{b_l};\frac{c_1}{d_1},...,\frac{c_l}{d_l}):\mathbb{B}
\rightarrow gpd$ for $a_1,...,a_k,b_1,...,b_k,\\
c_1,...,c_l,d_1,...,d_l \in \mathbb{N}^{+}.$

\begin{defi}
For $x \in \mathbb{B}$ we set
$$H(\frac{a_1}{b_1},...,\frac{a_k}{b_k};\frac{c_1}{d_1},...,\frac{c_l}{d_l})(x)=
\prod_{i=1}^{k}([a_i])_{|x|,[b_i]}\overline{\mathbb{Z}}_{b_i}^{x}
\prod_{j=1}^{l}{[d_j]}^{x}\overline{\mathbb{Z}}_{c_j,|x|,d_j}.$$
\end{defi}

\begin{thm}
$|H(\frac{a_1}{b_1},...,\frac{a_k}{b_l};\frac{c_1}{d_1},...,\frac{c_l}{d_l})|=
h(\frac{a_1}{b_1},...,\frac{a_k}{b_l};\frac{c_1}{d_1},...,\frac{c_l}{d_l})$.
\end{thm}

\begin{proof}
It follows from Lemma \ref{gps}, Lemma \ref{ya} and the fact that
$$H(\frac{a_1}{b_1},...,\frac{a_k}{b_l};\frac{b_1}{c_1},...,\frac{b_l}{c_l})=
\times_{i=1}^{k}H(\frac{a_i}{b_i},\emptyset)
\times_{j=1}^{l}H(\emptyset;\frac{c}{d}).$$
Since then $$|\times_{i=1}^{k}H(\frac{a_i}{b_i},\emptyset)
\times_{j=1}^{l}H(\emptyset;\frac{c}{d})|=\times_{i=1}^{k}h(\frac{a_i}{b_i},\emptyset)
\times_{j=1}^{l}h(\emptyset;\frac{c_j}{d_j})
=h(\frac{a_1}{b_1},...,\frac{a_k}{b_l};\frac{b_1}{c_1},...,\frac{b_l}{c_l}).$$
\end{proof}

According to Definition \ref{pf} the groupoid
$H(\frac{a_1}{b_1},...,\frac{a_k}{b_l};\frac{b_1}{c_1},...,\frac{b_l}{c_l})(x)$
may be described explicitly as the groupoid whose objects are
triples $(I,f,g)$ such that: 1) $I =(I_1,...,I_k)$ and $I_i
\subseteq [|x|-1].$ 2) $f= (f_1,...,f_k)$ and $f_i : \{0,...,|x|-1
\}
\rightarrow [a_i] \sqcup ([b_i]\times \mathbb{N})$ is such that
$f(i) \in [a_i]$ if $i \notin I$, and if $i \in I$ then $f(i) \in
[b_i]\times \mathbb{N}$ and  $1 \leq \pi_{\mathbb{N}}(f(i))\leq
i$. 3) $g=(g_1,...,g_l)$ where $g_j : x \rightarrow [d_j]$ for
$1\leq j \leq l.$ Morphisms are described as follows:
$$H(\frac{a_1}{b_1},...,\frac{a_k}{b_l};\frac{c_1}{d_1},...,\frac{c_l}{d_l})(x)((I,f,g),(J,p,q))=
\prod_{i=1}^{k}\overline{\mathbb{Z}}_{b_i}^{|x|}\prod_{j=1}^{l}\overline{\mathbb{Z}}_{c_j , |x|,
d_j},$$ if $I=J,$ $f_{i}|_{I_i} = p_{i}|_{I_i}$, and $g_j =q_j$;
otherwise,
$$H(\frac{a_1}{b_1},...,\frac{a_k}{b_l};\frac{c_1}{d_1},...,\frac{c_l}{d_l})(x)((I,f,g),(J,p,q))=
\emptyset.$$

\smallskip

Another combinatorial interpretation for
$h(\frac{a_1}{b_1},...,\frac{a_k}{b_l};\frac{b_1}{c_1},...,\frac{b_l}{c_l})$
is provided by the functor sending a finite set $x$ into

$$\prod_{i=1}^{k}(\overline{\mathbb{Z}}_{b_i}^{\oplus a_i})_{|x|}
\prod_{j=1}^{l}{[d_j]}^{x}\overline{\mathbb{Z}}_{c_j,|x|,d_j}.$$

 We close this paper by stating three
open problems. In this work we have provided a combinatorial
interpretation for hypergeometric functions, thus the following
problem seems natural.

\begin{problem}
Find a combinatorial interpretation for the $q$-hypergeometric
functions.
\end{problem}

In \cite{BD} we provide a combinatorial interpretation for the
Bernoulli and Euler numbers. Our interpretation of the Bernoulli
numbers seems to be related to the results of \cite{leechae}.

\begin{problem}
Find a combinatorial interpretation for the $q$-Bernoulli numbers
from \cite{Kim2}.
\end{problem}

\begin{problem}
Find a combinatorial interpretation for the $q$-Euler numbers.
\end{problem}

The interested reader will find more information on the various
$q$-analogues of the Bernoulli and Euler numbers in references
\cite{Kim}, \cite{Kim2}.

\subsection* {Acknowledgment}
This paper is dedicated to a true teacher and friend Professor
Crisitina Betz.

\end{document}